\documentclass[11pt,a4page]{article}
\usepackage{amssymb}
\usepackage{amsmath}
\usepackage{amsfonts}

\usepackage{amsthm}
\usepackage{bookmark}
\usepackage{cite}

\usepackage{changes}
\usepackage{color}
\usepackage{verbatim}

\date{}

\newtheorem{Theorem}{Theorem}[section]
\newtheorem{Lemma}{Lemma}[section]

\numberwithin{equation}{section} \theoremstyle{plain}

\def\R{{\textbf{R}}}

\def\f{\frac}
\def\R{\mathbb R}

\def\v{\mathbf{v}}

\def\al{\alpha}

\def\v{\varphi}

\title {Minimizers of $L^2$-Subcritical Inhomogeneous Variational Problems with A Spatially Decaying Nonlinearity}
\date{\today}
\author{ Yongshuai Gao ,  Yujin Guo\ and \   Shuang Wu \thanks{
		swu@mails.ccnu.edu.cn} \\
 {\small
          School of Mathematics and Statistics,
  Central China Normal University,}  \\
    { \small       Wuhan 430079, P. R. China}
       }
\author{ Yongshuai Gao\thanks{Email: ysgao@mails.ccnu.edu.cn.}, Yujin Guo\thanks{Email: yguo@ccnu.edu.cn. Y. J. Guo is partially supported by NSFC under Grant 11931012.
	} \ and \   Shuang Wu\thanks{Email: swu@mails.ccnu.edu.cn.}  \\
\small \it	School of Mathematics and Statistics,\\
\small \it  Hubei key Laboratory of Mathematical Sciences,\\
\small \it Central China Normal University,\\
\small \it Wuhan 430079, People's Republic of China\\}

\begin{document}
\baselineskip= 15pt
 \maketitle


\maketitle
\begin{abstract}
We study the minimizers of $L^2$-subcritical inhomogeneous variational problems with spatially decaying nonlinear terms, which contain $x = 0$ as a singular point. The limit concentration behavior of minimizers is proved as $M\to\infty$ by establishing the refined analysis of the spatially decaying nonlinear term.\\


\end{abstract}

\bigskip

\section{Introduction}
In this paper, we consider the minimizers of the following $L^2$-subcritical constraint inhomogeneous variational problem
\begin{equation}\label{A1problem}
I(M):=\underset{\{u\in \mathcal{H},\|u\|_2^2=1\}}{\inf}E_M(u),\ M>0,
\end{equation}
where the Gross-Pitaevskii (GP) energy functional $E_M(u)$ contains a spatially decaying nonlinearity and is defined by
\begin{equation}\label{A2functional}
E_M(u):=\int_{\R^N}\big(|\nabla u|^2+V(x)|u|^2\big)dx-\f{2M^{\f{p-1}{2}}}{p+1}\int_{\R^N}\f{|u|^{p+1}}{|x|^b}dx, \ N\ge 1,
\end{equation}
and the space $\mathcal{H}$ is defined as
\begin{equation*}
\mathcal{H}:=\Big\{u(x)\in H^1(\R^N):\ \int_{\R^N}V(x)|u(x)|^2<\infty\Big\}
\end{equation*}
with the associated norm $\|u\|_\mathcal{H}=\Big\{\displaystyle{\int_{\R^N}\Big(|\nabla u(x)|^2+|u(x)|^2+V(x)|u(x)|^2\Big)dx}\Big\}^\f{1}{2}$.
Here positive constants $b>0$ and $p>0$ of (\ref{A1problem}) satisfy
\begin{equation}\label{1:b}
0<b<\min\{2,N\},\ \ 1<p<1+\f{4-2b}{N},\ \ \mbox{where} \ \ N\ge 1,
\end{equation}
so that $E_M(u)$ admits $x = 0$ as a singular point in its nonlinear term. We always assume that the trapping potential $V(x)\geq0$ satisfies
\begin{enumerate}
\item [($V$).]  $V(x)\in L^\infty_{loc}(\R^N)\cap C_{loc}^\alpha(\R^N)$ with $\alpha\in(0,1)$, $\{x\in \R^N: V(x)=0\}=\{0\}$ and $\lim\limits_{|x|\to\infty}V(x)= \infty$.
\end{enumerate}
The variational problem (\ref{A1problem}) arises in various physical contexts, including the propagation of a laser beam in the optical fiber, Bose-Einstein condensates (BECs), and nonlinear optics (cf. \cite{A,B,LT}), where the constant $M>0$ often represents the attractive interaction strength, and $V(x)\ge 0$ denotes an external potential. The variational problem (\ref{A1problem}) and its associated elliptic equation have attracted a lot of attentions over the past few years, due to the appearance of the singular point $x=0$ in the nonlinear term, see \cite{AD,CG,dF,D,F,GS,S2} and the references therein.

When $b=0$, (\ref{A1problem}) is a homogeneous constraint variational problem, for which there are many existing results (\ref{A1problem}) (cf. \cite{C,FR,GLW,GSR,GWZZ,GZZ2,L1,L2,M,ZJ1}), including the existence and nonexistence of minimizers, and their quantitative properties of all kinds. More precisely, when $p>1+\f{4}{N}$, one can use the energy estimates to obtain the nonexistence of minimizers for (\ref{A1problem}) with $b=0$ as soon as $M>0$ (cf. \cite{C,CL}), which is essentially in the $L^2$-supercritical case. However, if $p=1+\f{4}{N}$, then (\ref{A1problem}) with $b=0$ reduces to the $L^2$-critical case, which was addressed widely by the second author and his collaborators, see \cite{GLW,GSR,GWZZ,GZZ2} and the references therein. As for the case where $1<p<1+\f{4}{N}$, (\ref{A1problem}) with $b=0$ is in the $L^2$-subcritical case and admits generally minimizers for all $M\in(0,\infty)$. In this case, the uniqueness, symmetry breaking and concentration behavior of minimizers were investigated recently as $M\to\infty$, see \cite{GZZ1,L1,M} and the references therein.

When $b\neq0$, the variational problem (\ref{A1problem}) contains the inhomogeneous nonlinear term $m(x)|u|^{p+1}$, where $m(x)=\f{1}{|x|^b}$ admits $x=0$ as a singular point. We remark that the inhomogeneous $L^2$-constraint variational problems were analyzed recently in \cite{DGL1,DGL2,M} and the references therein. However, as far as we know, the above mentioned works handle mainly with the inhomogeneous nonlinear term $m(x)|u|^{p+1}$
where $m(x)$ satisfies $m(x)\in L^\infty(\R^N)$ without any singular point. On the other hand, Ardila and Dinh obtained recently in \cite{AD} the existence of minimizers and the stability of the standing waves, for which they studied the associated constraint variational problem (\ref{A1problem}), in the $L^2$-subcritical case where the harmonic potential satisfies $V(x)=\gamma^2|x|^2 (\gamma>0)$, $b>0$ and $p>0$ satisfy (\ref{1:b}).

Under the assumptions $(V)$ and (\ref{1:b}), we comment that it is standard to obtain the existence of minimizers for $I(M)$ for all $M > 0$, see \cite [Theorem 1.8] {AD} and the related argument. Motivated by  above mentioned works, in this paper we mainly study the limit behavior of minimizers $u_M$ for $I(M)$ as $M\to\infty$, and the main purpose of
this paper is to investigate the impact of the singular point $x=0$ on the behavior of
$u_M$ as $M\to\infty$.

We now assume that $u_M$ is a minimizer of $I(M)$ for any $M > 0$. It then follows
from the variational theory that $u_M$ satisfies the following Euler-Lagrange equation
\begin{equation}\label{A8ELE}
  -\Delta u_M+V(x)u_{M}-M^{\f{p-1}{2}}\frac{u_M^p}{|x|^{b}}=\mu_{M}u_{M}\ \ \text{in}\ \ \R^N,
\end{equation}
where $\mu_{M}\in\R^N$ is a suitable Lagrange multiplier associated to $u_{M}$. By the form of the energy functional $E_M(\cdot)$, one can obtain from \cite [Theorem 6.17] {LL} that $E_M(u)=E_M(|u|)$ holds for any $u\in \mathcal{H}$, which implies that $|u_M|$ is also a minimizer of $I(M)$. By the strong maximum principle, one can further derive from (\ref{A8ELE}) that $|u_M| > 0$ holds in $\R^N$. Therefore, $u_M$ must be either positive or negative. Without loss of generality, in the following we only consider positive minimizers $u_M >0$ of $I(M)$.

Under the assumption (\ref{1:b}), we next recall the following sharp Gagliardo-Nirenberg (GN) inequality (cf. \cite [Theorem 1.2] {F}):
\begin{equation}\label{A3GNineq}
\int_{\R^N}\frac{|u|^{p+1}}{|x|^b}dx\leq C_{GN}^{-1}\|\nabla u\|_{2}^{\f{N(p-1)}{2}+b}\|u\|_{2}^{p+1-\f{N(p-1)}{2}-b},\ \ u\in H^1(\R^N),
\end{equation}
where $C_{GN}>0$ satisfies
\begin{equation}\label{A4GNconstant}
C_{GN}=\Big(\f{N(p-1)+2b}{2(p+1)-N(p-1)-2b}\Big)^{\f{N(p-1)+2b}{4}}\f{2(p+1)-N(p-1)-2b}{2(p+1)}\|w\|_2^{p-1},
\end{equation}
and $w$ is the unique positive radially symmetric solution (cf. \cite{BL,GS,GNN,LN,Y}) of
\begin{equation}\label{A5wdeequation}
-\Delta w+w-\frac{w^p}{|x|^{b}}=0\,\ \text{in}\, \ \R^N,\ \ w\in H^1(\R^N).
\end{equation}
The equality in (\ref{A3GNineq}) is achieved at $u=w$. Moreover, $w$ satisfies the following Pohozaev identity
\begin{equation}\label{A7wdeidentity}
\|\nabla w\|_{2}^2=\f{N(p-1)+2b}{2(p+1)}\int_{\R^N}\f{|w|^{p+1}}{|x|^{b}}dx=\f{N(p-1)+2b}{2(p+1)-N(p-1)-2b}\|w\|_{2}^2.
\end{equation}
Recall also from \cite [Theorem 2.2] {GS} that there exist  positive constants $\delta>0$ and $C>0$ such that $w(x)$ satisfies
\begin{equation}\label{A6wdedecay}
w(x),\ |\nabla w(x)|\leq Ce^{-\delta|x|}\ \ \text{as}\, \ |x|\to\infty.
\end{equation}
All above properties of $w$ are often used in the refined analysis of minimizers for $I(M)$
as $M\to\infty$.

Using above notations, the main result of the present paper can be stated as
the following theorem.


\begin{Theorem}\label{limit.theorem}
Under the assumptions $(V)$ and (\ref{1:b}), let $u_k$ be a positive minimizer of $I(M_k)$, where $M_k \to \infty$ as $k\to \infty$. Then there exists a subsequence, still denoted by $\{u_k\}$, of $\{u_k\}$ such that $u_{k}$ satisfies
\begin{equation}\label{A9limit.behavior}
  w_k(x):=\epsilon_k^{\f{N}{2}}u_k(\epsilon_kx)\to \f{w(x)}{\sqrt{a^*}}\ \ \text{uniformly in}\ \  L^\infty(\R^N) \,\ \text{as}\, \  k\to \infty,
\end{equation}
where $\epsilon_k:=\Big({\frac{M_k}{a^*}}\Big)^{-\frac{p-1}{4-N(p-1)-2b}}>0$, $a^*:=\|w\|_2^2>0$, and $w>0$ is the unique positive solution of \eqref{A5wdeequation}.
Moreover, $u_k$ decays exponentially in the sense that for sufficiently large $k>0,$
\begin{equation}\label{eqSA2}
w_k(x)\leq Ce^{-\sqrt{\theta}|x|} \ \ \text{and}\ \ |\nabla w_k(x)|\leq Ce^{-\theta|x|}\ \ \text{as} \ \ |x|\to\infty,
\end{equation}
where $0<\theta<1$ and $C>0$ are independent of $k>0$.
\end{Theorem}

The proof of Theorem \ref{limit.theorem} shows essentially that as $M_k\to\infty$, $u_k$ prefers to concentrate near the singular point $x=0$ of $I(M_k)$, instead of a minimum point for $V(x)$. The proof of Theorem \ref{limit.theorem} depends on the refined estimate of $\int_{\R^N} V (x)|u_k|^2dx$ as $k\to\infty$, for which we shall consider the following constraint variational problem without the
trap:
\begin{equation}\label{A12noVproblem}
  \tilde{I}(M):=\underset{\{u\in H^1(\R^N),\,\|u\|_2^2=1\}}{\inf}\tilde{E}_M(u),\,\ N\geq1,
\end{equation}
where $\tilde{E}_M(u)$ is defined by
\begin{equation}\label{A13noVfunctional}
 \tilde{E}_M(u):=\int_{\R^N}|\nabla u|^2dx-\f{2M^{\f{p-1}{2}}}{p+1}\int_{\R^N}\f{|u|^{p+1}}{|x|^b}dx.
\end{equation}
By deriving the energy
estimates between $\tilde{I}(M_k)$ and $I(M_k)$ as $M_k\to\infty$, we shall verify that $I(M_k)-\tilde{I}(M_k)\to0$ as $M_k\to\infty$, which further implies that $\int_{\R^N} V (x)|u_k|^2dx\to0$ as $M_k\to\infty$. Furthermore, the $L^\infty$-uniform convergence of (\ref{A9limit.behavior}), which is established by analyzing delicately the singular nonlinear term of $I(M_k)$, seems crucial in the further refined investigations on the minimizers of $I(M_k)$ as $M_k\to\infty$.

This paper is organized as follows. Section 2 is devoted to the refined energy estimates of $I(M)$ as $M\to\infty$, based on which we shall complete in Section 3 the proof of Theorem \ref{limit.theorem} on the limit behavior of minimizers for $I(M)$ as $M\rightarrow\infty$.

\section{Energy estimates of $I(M)$}

This section is devoted to establishing the energy estimates of $I(M)$ as $M\to\infty$ by analyzing the energy of $\tilde{I}(M)$ defined in (\ref{A12noVproblem}). Employing the concentration-compactness principle, one can deduce that $\tilde{I}(M)$ admits minimizers for any $M\in(0,\infty)$, see, e.g., \cite{C,L1,L2}. Moreover, without loss of generality, we may consider positive minimizers of $\tilde{I}(M)$ defined in (\ref{A12noVproblem}). We start with the following energy estimates of $\tilde{I}(M)$.

\begin{Lemma}\label{LemmanoVenergy.esti}
Under the assumption (\ref{1:b}), assume that $V(x)$ satisfies $(V)$, and let $\tilde {u}_M$ be a positive minimizer of $\tilde{I}(M)$. Then for any $M>0,$
\begin{equation}\label{noVenergy.esti}
\tilde{I}(M)=-\lambda_0\Big (\f {M}{a^*}\Big)^{\f {2(p-1)}{4-N(p-1)-2b}},
\end{equation}
       and
\begin{equation}\label{noVminimizer}
\tilde{u}_M(x)=\f{1}{\sqrt{a^*}}\tilde{\alpha}^{\f{N}{2}}_Mw(\tilde{\alpha}_Mx),
\end{equation}
where $\tilde{\alpha}_M:=\big(\f {M}{a^*}\big)^{\f{p-1}{4-N(p-1)-2b}}>0$, $\lambda_0:=-\f {N(p-1)+2b-4}{2(p+1)-N(p-1)-2b}>0$ and $a^*:=\|w\|_2^2$. Here $w>0$ is the unique positive solution of the equation (\ref{A5wdeequation}).
\end{Lemma}

\noindent{\bf Proof.}
	   Assume that $\tilde{u}_M$ is a positive minimizer of $\tilde{I}(M)$ and  $\tilde{u}_1$ is a positive minimizer of $\tilde{I}(1)$. We claim that  for any $M>0,$
\begin{equation}\label{eqS2.3}
\tilde{I}(M)=M^{\f{2(p-1)}{4-N(p-1)-2b}}\tilde{I}(1) \text{ and } \tilde{u}_M(x)=\alpha _M^\f N2 \tilde{u}_1(\alpha _Mx),
\end{equation}
where $\alpha_M:=M^{\f{p-1}{4-N(p-1)-2b}}>0$. Indeed, setting $\tilde{w}_1(x):=\al _M^{-\f{N}{2}}\tilde{u}_M(\al ^{-1}_Mx)$, one can deduce from (\ref{A12noVproblem}) that
\begin{align}\label{G1}
\tilde{I}(M)=\tilde{E}_M(\tilde{u}_M)
&=\int_{\R^N}|\nabla\tilde{u}_M|^2dx-\f{2M^{\f{p-1}{2}}}{p+1}\int_{\R^N}\f{|\tilde{u}_M|^{p+1}}{|x|^b}dx\notag\\
&=\al _M^2\int _{\R^N}|\nabla \tilde{w}_1|^2dx-\f{2M^{\f{p-1}{2}}}{p+1}\cdot\al _M^{\f{N(p+1)}{2}}\cdot \al _M^{-N}\cdot \al _M^b\int_ {\R^N}\f{|\tilde{w}_1|^{p+1}}{|x|^b}dx\notag\\
&=M^{\f{2(p-1)}{4-N(p-1)-2b}} \Big[ \int_{\R^N}|\nabla\tilde{w}_1|^2dx-\f{2}{p+1}\int_ {\R^N}\f{|\tilde{w}_1|^{p+1}}{|x|^b}dx\Big]\notag\\
&\geq M^{\f{2(p-1)}{4-N(p-1)-2b}}\tilde{I}(1).
\end{align}
Similarly, setting $\tilde{w}_M(x):=\al _M^{\f{N}{2}}\tilde{u}_1(\al _Mx)$ as a test function of $\tilde{I}(M)$, one can get that
\begin{equation}\label{G2}
  \tilde{I}(M) \leq \tilde{E}_M(\tilde{w}_M) =M^{\f{2(p-1)}{4-N(p-1)-2b}} \tilde{I}(1).
\end{equation}
Following \eqref{G1} and \eqref{G2}, we conclude that the first equality of (\ref{eqS2.3}) holds. Furthermore, one can check that $\tilde{w}_1$ is a minimizer of $\tilde{I}(1)$ and $\tilde{w}_M$ is a minimizer of $\tilde{I}(M)$. This proves the second equality of (\ref{eqS2.3}). Therefore, the claim (\ref{eqS2.3}) holds true.

We next prove that for any $M>0$,
\begin{equation}\label{eqS2.4}
\tilde{I}(1)=-\lambda_0(a^*)^{-\f {2(p-1)}{4-N(p-1)-2b}},\,\ \text{where} \,\ \lambda_0:=-\f {N(p-1)+2b-4}{2(p+1)-N(p-1)-2b}>0,
\end{equation}
and
\begin{equation}\label{eqS2.5}
\tilde{u}_1(x)={(a^*)}^{-\f{2-b}{4-N(p-1)-2b}}w\big((a^*)^{-\f{p-1}{4-N(p-1)-2b}}x\big).
\end{equation}
Consider  a test function $0<\tilde{v}_0\in H^1(\R^N)$ satisfying $||\tilde{v}_0||_2^2=1$. Set $\tilde{v}_{\epsilon}(x):=\epsilon^{\f{N}{2}}\tilde{v}_0(\epsilon x)$, where $\epsilon>0$ is small enough. One can get that for sufficiently small $\epsilon>0$,
\begin{align}\label{eqS2.6}
\tilde{I}(1)\leq \tilde{E}_1(\tilde{v}_{\epsilon})&=\int_{\R^N}|\nabla\tilde{v}_{\epsilon}|^2dx-\f{2}{p+1}
\int_{\R^N}\f{|\tilde{v}_\epsilon|^{p+1}}{|x|^b}dx\notag	
\\&=\epsilon^2\int_{\R^N}|\nabla\tilde{v}_0|^2dx-\f{2\epsilon^{\f{N(p-1)}{2}+b}}{p+1}
\int_{\R^N}\f{|\tilde{v}_0|^{p+1}}{|x|^b}dx<0,
\end{align}
due to the assumption (\ref{1:b}).
Let $\tilde{u}_1>0$ be a positive minimizer of $\tilde{I}(1)$. Then $\tilde{u}_1$ satisfies the following Euler-Lagrange equation
\begin{equation}\label{eqS2.7}
   -\Delta\tilde{u}_1(x)=\tilde{\mu}_1\tilde{u}_1(x)+\f{\tilde{u}_1^p(x)}{|x|^b}\ \ \text{in} \,\ \R^N,
\end{equation}
where $\tilde{\mu}_1\in\R$ is the Lagrangian multiplier associated to $\tilde{u}_1$.  Applying (\ref{eqS2.6}) and (\ref{eqS2.7}), we get that
\begin{equation}\label{eqS2.8}
\arraycolsep=1.5pt\begin{array}{lll}
\tilde{\mu}_1&=\displaystyle\int_{\R^N}|\nabla\tilde{u}_1|^2dx-\int_{\R^N}\f{|\tilde{u}_1|^{p+1}}{|x|^b}dx\\[3mm]
&=\displaystyle\tilde{I}(1)-\f{p-1}{p+1}\int_\R\f{|\tilde{u}_1|^{p+1}}{|x|^b}dx<0.
\end{array}
\end{equation}
Since $w>0$ is the unique positive solution of (\ref{A5wdeequation}), one can conclude from (\ref{A5wdeequation}) and (\ref{eqS2.7}) that
\begin{equation*}
 \tilde{u}_1(x)=(-\tilde{\mu}_1)^{\f{2-b}{2(p-1)}}w\big(({-\tilde{\mu}_1})^{\f{1}{2}}x\big),
\end{equation*}
where $\tilde{\mu}_1<0$ holds by (\ref{eqS2.8}). Moreover, since
$$1=\|\tilde{u}_1\|_2^2
=(-\tilde{\mu}_1)^{\f{4-2b-N(p-1)}{2(p-1)}}\|w\|_2^2=(-\tilde{\mu}_1)^{\f{4-2b-N(p-1)}{2(p-1)}}a^*,$$
one can derive that
\begin{equation*}
 \tilde{\mu}_1=-(a^*)^{\f{2(1-p)}{4-2b-N(p-1)}}<0 \ \ \text{ and } \ \ \tilde{u}_1(x)=(a^*)^{\f{b-2}{4-2b-N(p-1)}}w\Big((a^*)^{\f{1-p}{4-2b-N(p-1)}}x\Big)>0.
\end{equation*}
Hence, (\ref{eqS2.5}) is proved. On the other hand, substituting (\ref{eqS2.5}) and (\ref{A7wdeidentity}) into (\ref{A13noVfunctional}), we get that
\begin{align*}
\tilde{I}(1)=\tilde{E}_1(\tilde{u}_1)&=\int_{\R^N}|\nabla\tilde{u}_1|^2dx-\f{2}{p+1}\int_{\R^N}\f{|\tilde{u}_1|^{p+1}}{|x|^b}dx\\
&=(a^*)^{\f{2(b-2)}{4-2b-N(p-1)}}\cdot (a^*)^{\f{2(1-p)}{4-2b-N(p-1)}}\cdot (a^*)^{\f{N(p-1)}{4-2b-N(p-1)}}\int_{\R^N}|\nabla w|^2dx\\
&\quad -(a^*)^{\f{(b-2)(p+1)}{4-2b-N(p-1)}}\cdot(a^*)^{\f{N(p-1)}{4-2b-N(p-1)}}\cdot(a^*)^{\f{(1-p)b}{4-2b-N(p-1)}}\cdot \f{2}{p+1}\int_{\R^N}\f{|w|^{p+1}}{|x|^b}dx\\
&=(a^*)^{\f{2(b-2)+(p-1)(N-2)}{4-2b-N(p-1)}}\cdot \f{N(p-1)+2b}{2(p+1)-N(p-1)-2b}\cdot a^*\\
&\quad-(a^*)^{\f{2(b-2)+(p-1)(N-2)}{4-2b-N(p-1)}}\cdot \f{2}{p+1}\cdot \f{2(p+1)}{2(p+1)-N(p-1)-2b}\cdot a^*\\
&=(a^*)^{-\f{2(p-1)}{4-2b-N(p-1)}}\f{N(p-1)+2b-4}{2(p+1)-N(p-1)-2b}\\
&=-\lambda_0(a^*)^{-\f{2(p-1)}{4-2b-N(p-1)}},
\end{align*}
which thus implies that (\ref{eqS2.4}) holds.

We finally conclude from (\ref{eqS2.3})--(\ref{eqS2.5}) that for any $M>0$,
\begin{equation*}
   \tilde{I}(M)=M^{\f{2(p-1)}{4-N(p-1)-2b}}\cdot (a^*)^{-\f{2(p-1)}{4-N(p-1)-2b}}\cdot (-\lambda_0)=-\lambda_0\Big(\f{M}{a^*}\Big)^{\f{2(p-1)}{4-N(p-1)-2b}},
\end{equation*}
and
\begin{align*}
 \tilde{u}_M(x)=\al_M^{\f{N}{2}}\tilde{u}_1(\al_Mx)&=\al_M^{\f{N}{2}}(a^*)^{-\f{2-b}{4-N(p-1)-2b}}
 w\Big((a^*)^{-\f{p-1}{4-N(p-1)-2b}}\cdot \al_Mx\Big)\\
 &=M^{\f{N(p-1)}{2[4-N(p-1)-2b]}}\cdot (a^*)^{-\f{2-b}{4-N(p-1)-2b}}\cdot w\Big((\f{M}{a^*})^{\f{p-1}{4-N(p-1)-2b}} x\Big)\\
 &=\f{1}{\sqrt{a^*}}\tilde{\alpha}_M^{\f{N}{2}}w(\tilde{\alpha}_Mx).
\end{align*}
Therefore, the proof of Lemma \ref{LemmanoVenergy.esti} is completed.
 \qed

Applying Lemma \ref{LemmanoVenergy.esti}, we now establish the energy estimates of $I(M)$.

\begin{Lemma}\label{Lemmaenergy.esti}
Under the assumption (\ref{1:b}), assume that $V(x)$ satisfies $(V)$. Then we have
\begin{equation}\label{energy.esti}
\lim\limits_{M\to \infty}\f{I(M)}{(\f{M}{a^*})^{\f{2(p-1)}{4-N(p-1)-2b}}}=-\lambda_0,
\end{equation}
where $a^*:=\|w\|_2^2$, $\lambda_0:=-\f {N(p-1)+2b-4}{2(p+1)-N(p-1)-2b}>0$, and $w>0$ is the unique positive solution of the equation (\ref{A5wdeequation}).
\end{Lemma}

\noindent{\bf Proof.}
We first establish the lower bound of $I(M)$ as $M\to\infty$. Let $u_M>0$ be a positive minimizer of $I(M)$. Under the assumption $(V)$, we get from (\ref{A12noVproblem}) that
\begin{align}\label{energy.lower.esti}
I(M)&=\int_{\R^N}|\nabla u_M|^2dx+\int_{\R^N}V(x)u_M^2(x)dx-\f{2M^{\f{p-1}{2}}}{p+1}\int_{\R^N}\f{|u_M|^{p+1}}{|x|^b}dx\notag\\
&\geq \int_{\R^N}|\nabla u_M|^2dx-\f{2M^{\f{p-1}{2}}}{p+1}\int_{\R^N}\f{|u_M|^{p+1}}{|x|^b}dx\\
&\geq \widetilde{I}(M)=-\lambda_0\Big(\f{M}{a^*}\Big)^{\f{2(p-1)}{4-(p-1)-2b}} \ \ \mbox{as}\ \ M\to \infty,\notag
\end{align}
where $\lambda_0>0$ is as in Lemma \ref{LemmanoVenergy.esti}. This thus implies  the lower bound of $I(M)$ as $M\to\infty$.

We next estimate the upper bound of $I(M)$ as $M\to\infty$. Define
\begin{equation}\label{eqS2.11}
u_{\tau}(x):=\f{A_{\tau}\tau^{\f{N}{2}}}{\|w\|_2}w(\tau x)\v (x),\,\ \tau>0,
\end{equation}
where $0\le \v(x)\in C^{\infty}(\R^N)$ is a cut-off function satisfying
$$ \v (x)=
    \begin{cases}
    1,\quad\quad |x|\leq 1;\\ 0,\quad\quad |x|\geq 2,
    \end{cases}$$
$w>0$ is the unique positive solution of (\ref{A5wdeequation}), and $A_\tau>0$ is a suitable constant such that $\|u_{\tau}(x)\|_2^2$=1. Applying the exponential decay of $w$ in (\ref{A6wdedecay}), one can check that as $\tau\to\infty$,
\begin{align}\label{eqS2.12}
   1\leq A^2_{\tau}=\f{\|w\|^2_2}{\int_{\R^N}w^2(x)\v ^2(\f{x}{\tau})dx}
   \leq 1+\f{\int_{B^c_\tau}w^2(x)dx}{\int_{B_\tau}w^2(x)dx}
   \leq 1+Ce^{-2\delta \tau},
\end{align}
where $C>0$ is independent of $\tau>0$.
Substituting (\ref{eqS2.11}) into \eqref{A2functional} and applying the exponential decay of $w$ in \eqref{A6wdedecay} and the identity (\ref{A7wdeidentity}), direct calculations yield that as $\tau\to\infty$,
\begin{align}\label{eqS2.13}
    \int_{\R^N}|\nabla u_\tau|^2dx
    &=\f{A^2_{\tau}}{\|w\|^2_2}\int_{\R^N}\Big|\v \big(\f{x}{\tau}\big)\tau\nabla w(x)+w(x)\nabla \v \big(\f{x}{\tau}\big)\Big|^2dx\nonumber\\
    &=\f{A^2_{\tau}}{\|w\|^2_2}\int_{\R^N}\bigg[\tau^2\v ^2\big(\f{x}{\tau}\big)\big|\nabla w(x)\big|^2
    +w^2(x)\big|\nabla \v \big(\f{x}{\tau}\big)\big|^2
    \\
    &\qquad \qquad\qquad\ +2\tau\nabla w(x)\v\big(\f{x}{\tau}\big)
    \nabla\v \big(\f{x}{\tau}\big)w(x)\bigg] dx\nonumber\\
    &\leq \big(1+Ce^{-2\delta\tau}\big)\f{N(p-1)+2b}{2(p+1)-N(p-1)-2b}\tau^2,\nonumber
\end{align}
\begin{align}\label{eqS2.14}
    \int_{\R^N}V(x)u^2_{\tau}(x)dx
    &=\f{A^2_{\tau}}{\|w\|^2_2}\int_{\R^N}V\big(\f{x}{\tau}\big)w^2(x)\v^2\big(\f{x}{\tau}\big)dx\nonumber\\
    &=\f{A^2_{\tau}}{\|w\|^2_2}\int_{B_{2\tau}(0)}V\big(\f{x}{\tau}\big)w^2(x)\v^2\big(\f{x}{\tau}\big)dx\\
    &\leq \f{A^2_{\tau}}{\|w\|^2_2}\Big[\int_{0\leq|x|\leq\sqrt{\tau}}V(\f{x}{\tau})w^2(x)dx+Ce^{-\delta \sqrt{\tau}}\Big],\nonumber
\end{align}
and
\begin{align}\label{eqS2.15}
    \int_{\R^N}\f{|u_{\tau}|^{p+1}}{|x|^b}dx
    &=\f{A^{p+1}_\tau \tau^{\f{N}{2}(p+1)}}{\|w\|^{p+1}_2}
    \int_{\R^N}\f{|w(\tau x)\v(x)|^{p+1}}{|x|^b}dx\nonumber\\
    &=\f{A^{p+1}_\tau \tau^{\f{N}{2}(p-1)+b}}{\|w\|^{p+1}_2}
    \int_{\R^N}\f{|w(x)\v(\f{x}{\tau})|^{p+1}}{|x|^b}dx\\
    &\geq \tau^{\f{N}{2}(p-1)+b}
    \f{2(p+1)}{2(p+1)-N(p-1)-2b}(a^*)^{\f{1-p}{2}}-Ce^{-\delta(p+1)\tau}\tau ^{\f{N(p-1)}{2}}.\nonumber
\end{align}
It then follows from (\ref{eqS2.13})--(\ref{eqS2.15}) that as $\tau\to\infty$,
\begin{align*}
I(M)\leq E_M(u_\tau)\le &\,\f{N(p-1)+2b}{2(p+1)-N(p-1)-2b}\tau^2
-\tau^{\f{N}{2}(p-1)+b}
    \f{4(\f{M}{a^*})^{\f{p-1}{2}}}{2(p+1)-N(p-1)-2b}\\
& +V(0)+Ce^{-\delta \sqrt{\tau}}.
\end{align*}
Setting $\tau=(\f{M}{a^*})^{\f{p-1}{4-N(p-1)-2b}}$ into the above estimate, it then gives that as $M\to \infty$,
\begin{equation}\label{energy.supper.esti}
I(M)\leq-\lambda_0\Big(\f{M}{a^*}\Big)^{\f{2(p-1)}{4-N(p-1)-2b}}+Ce^{-\delta \sqrt{\tau}}=-\lambda_0\Big(\f{M}{a^*}\Big)^{\f{2(p-1)}{4-N(p-1)-2b}}+o(1),
\end{equation}
where $\lambda_0:=-\f {N(p-1)+2b-4}{2(p+1)-N(p-1)-2b}>0$. Thus, (\ref{energy.esti}) follows from (\ref{energy.lower.esti}) and (\ref{energy.supper.esti}), which completes the proof of Lemma \ref{Lemmaenergy.esti}.
\qed

\section{Proof of Theorem 1.1}
In this section, we shall complete the proof of Theorem 1.1 on the limit behavior of minimizers for (\ref{A1problem}) by the blow-up analysis. We first establish the following lemma.

\begin{Lemma}\label{LemmaVto0}
Under the assumption (\ref{1:b}), assume that $V(x)$ satisfies $(V)$, and let $u_M$ be a positive minimizer of $I(M)$. Then we have
\begin{equation}\label{Vto0}
\int_{\R^N}V(x)u_M^2(x)dx\to 0\ \ \text{as}\, \ M\to\infty.
\end{equation}
\end{Lemma}

\noindent{\bf Proof.}
The key of proving \eqref{Vto0} is to verify that
\begin{equation}\label{eqS3.2}
I(M)-\tilde{I}(M)\to0 \ \ \text{as}\, \ M\to\infty.
\end{equation}
Indeed, if (\ref{eqS3.2}) holds, then one derive from (\ref{A12noVproblem}) and (\ref{A13noVfunctional}) that
\begin{equation}\label{eqS3.3}
\int_{\R^N}V(x)u_M^2(x)dx=I(M)-\tilde{E}_M(u_M)\leq I(M)-\tilde{I}(M)\, \to\, 0\ \ \text{as}\, \ M\to\infty,
\end{equation}
which thus implies that (\ref{Vto0}) holds.

We now prove (\ref{eqS3.2}). By Lemma \ref{LemmanoVenergy.esti}, we deduce from (\ref{energy.supper.esti}) that
\begin{equation}\label{eqS3.4}
I(M)\leq \tilde{I}(M)+o(1)\ \ \text{as}\, \ M\to\infty.
\end{equation}
On the other hand, it follows from (\ref{A12noVproblem}) and (\ref{A13noVfunctional}) that
\begin{equation}\label{eqS3.5}
I(M)-\tilde{I}(M)\geq E_M(u_M)-\tilde{E}_M(u_M)=\int_{{\R^N}}V(x)u_M^2(x)dx\geq 0\ \ \text{as}\, \ M\to\infty.
\end{equation}
Therefore, the estimate (\ref{eqS3.2}) now follows from (\ref{eqS3.4}) and (\ref{eqS3.5}), and we are done.
\qed

Motivated by \cite{GWZZ,GZZ2,WXF}, we next establish the following lemma.

\begin{Lemma}\label{lemma3.2}
Under the assumption (\ref{1:b}), assume that $V(x)$ satisfies $(V)$, and let $u_k$ be a positive minimizer of $I(M_k)$, where $M_k\to\infty$ as $k\to \infty$. Define
\begin{align}\label{eqS3.6}
w_k(x):=\epsilon_k^{\f N2}u_k\big(\epsilon_kx\big),
\end{align}
where $\epsilon_k:=(\f{M_k}{a^*})^{-\f{p-1}{4-N(p-1)-2b}}\to 0$ as $k\to\infty $.
Then there exists a subsequence, still denoted by $\{w_k\}$, of $\{w_k\}$ such that
\begin{equation}\label{eqS3.7}
w_k(x)\to \f{w(x)}{\sqrt{a^*}}\, \ \text{strongly in}\, \ H^1(\R^N)\, \ \text{as} \,\ k\to\infty,
\end{equation}
where $a^*:=\|w\|^2_2$ and $w>0$ is the unique positive solution of (\ref{A5wdeequation}).
\end{Lemma}

\noindent{\bf Proof.}
We first prove that there exist some positive constants $C_1$, $C_2$, $C'_1$, $C'_2$, which are independent of $k$, such that as $k\to\infty$,
\begin{equation}\label{eqS3.8}
0<C_1\leq \|\nabla w_k\|^2_2\leq C_2\ \text{ and }\ 0<C'_1\leq \int_{\R^N}\f{|{w}_k|^{p+1}}{|x|^b}dx \leq C'_2.
\end{equation}
Indeed, using Lemmas \ref{Lemmaenergy.esti} and \ref{LemmaVto0}, we deduce from \eqref{A2functional} and (\ref{eqS3.6}) that as $k\to \infty$,
\begin{equation}\label{eqS3.9}
\begin{split}
\epsilon_k^2I(M_k)
  &=\epsilon_k^2\Big(\int_{\R^N}|\nabla u_k|^2dx+\int_{\R^N}V(x)u_k^2dx-\f{2M_k^{\f{p-1}{2}}}{p+1}\int_{\R^N}\f{|u_k|^{p+1}}{|x|^b}dx\Big)\\
  &=\Big(\int_{\R^N}|\nabla{w}_k|^2dx
  -\f{2(a^*)^{\f{p-1}{2}}}{p+1}\int_{\R^N}\f{|{w}_k|^{p+1}}{|x|^b}dx
  +o(\epsilon^2_k )\Big)\to -\lambda_0<0,
\end{split}
\end{equation}
which implies that
\begin{equation}\label{eqS3.12}
 \f{2(a^*)^{\f{p-1}{2}}}{p+1} \f{\int_{\R^N}\f{|{w}_k|^{p+1}}{|x|^b}dx}{\int_{\R^N}|\nabla{w}_k|^2dx}\to 1\,\ \text{ as } \ k\to \infty.
\end{equation}
By contradiction, assume that $\|\nabla{w}_k\|_2^2\to \infty$ as $k\to \infty$. Define $\gamma_k^2:=\|\nabla{w}_k\|_2^2>0$ and $\upsilon_k (x)=\gamma_k^{-\f{N}{2}}{w}_k(\gamma_k^{-1}x)$, so that $\gamma_k^2\to \infty$ as $k\to \infty$. It then follows that $\|\upsilon_k\|^2_2=1$ and $\|\nabla\upsilon_k\|_2^2=1$ for all $k\ge 1$. Further, we deduce from the GN inequality (\ref{A3GNineq}) that
\begin{equation}\label{eqS3.10}
\int_{\R^N}\f{|\upsilon_k|^{p+1}}{|x|^b}dx\leq  C^{-1}_{GN}\|\nabla\upsilon_k\|^{\f{N(p-1)}{2}+b}_2\|\upsilon_k\|_{2}^{p+1-\f{N(p-1)}{2}-b}=C^{-1}_{GN},
\end{equation}
where $C_{GN}>0$ is given in (\ref{A4GNconstant}). Under the assumption (\ref{1:b}), since $\|v_k\|^2_2=\|\nabla v_k\|^2_2=1$ for all $k\geq 1$, it follows from (\ref{eqS3.10}) that
\begin{equation}\label{eqS3.11}
\begin{split}
 \f{2(a^*)^{\f{p-1}{2}}}{p+1}\f{\int_{\R^N}\f{|{w}_k|^{p+1}}{|x|^b}dx}{\int_{\R^N}|\nabla{w}_k|^2dx}
 &=\f{2(a^*)^{\f{p-1}{2}}}{p+1}\f{\int_{\R^N}\f{|\upsilon_k|^{p+1}}{|x|^b}dx}{\|\nabla\upsilon_k\|_2^2}\gamma_k^{-\f{4-N(p-1)-2b}{2}}\\
 &\le \f{2(a^*)^{\f{p-1}{2}}}{p+1}C^{-1}_{GN}\gamma_k^{-\f{4-N(p-1)-2b}{2}}\to 0\,\ \text{ as }\, \ k\to \infty,
\end{split}
\end{equation}
which however contradicts to (\ref{eqS3.12}). Hence, we conclude that $\|\nabla{w}_k\|_2^2\leq C_2$ holds uniformly as $k\to \infty$. Applying the GN inequality (\ref{A3GNineq}) and the fact that $\|w_k\|^2_2=1$, we deduce from above that $\int_{\R^N}\f{|{w}_k|^{p+1}}{|x|^b}dx \leq C'_2$ holds uniformly as $k\to \infty$. On the other hand, one can obtain from (\ref{eqS3.9}) that $\int_{\R^N}\f{|{w}_k|^{p+1}}{|x|^b}dx\geq C_1'$ holds uniformly as $k\to \infty$, together with (\ref{A3GNineq}), which then imply that $\|\nabla w_k\|^2_2\geq C_1$ holds uniformly as $k\to \infty$. We therefore conclude that (\ref{eqS3.8}) holds.

From (\ref{eqS3.6}) and (\ref{eqS3.8}), we deduce that $w_k$ is bounded uniformly in $H^1(\R^N)$, which implies that there exist a subsequence, still denoted by $\{w_k\}$, of $\{w_k\}$ and $0\le w_0\in H^1(\R^N)$ such that
 \begin{equation}\label{eqS3.13}
 w_k\rightharpoonup w_0 \geq 0\ \ \text{weakly in}\, \   H^1(\R^N)\, \ \text{as }\ k\to\infty.
 \end{equation}
We now prove that $w_0\not\equiv0$. Motivated by \cite{AD,GS}, we first claim that
 \begin{equation}\label{eqS3.18}
 \int_{\R^N}|x|^{-b}w_k^{p+1}dx\to\int_{\R^N}|x|^{-b}w_0^{p+1}dx\ \ \text{as}\ \ k\to\infty.
 \end{equation}
Actually, we have
\begin{align*}	
&\Big|\int_{\R^N}|x|^{-b}w_k^{p+1}dx-\int_{\R^N}|x|^{-b}w_0^{p+1}dx\Big|\\
\leq&\int_{\R^N}|x|^{-b}\Big|w_k^{p+1}-w_0^{p+1}\Big|dx\\
=&\int_{B_R}|x|^{-b}\Big|w_k^{p+1}-w_0^{p+1}\Big|dx
+\int_{\R^N\setminus
B_R}|x|^{-b}\Big|w_k^{p+1}-w_0^{p+1}\Big|dx:=A_k+B_k,
\end{align*}
where $R>0$ is arbitrary. Under the assumption (\ref{1:b}), we have $1<p<1+\f{4-2b}{N}$ and
$$ 1+\f{4-2b}{N}<1+ \f{4-2b}{N-2},\ \ \mbox{if}\, \ N\ge 3.$$
By H\"{o}lder inequality, we then have
\begin{align}	\label{eqS3.18A}
A_k=&\int_{B_R}|x|^{-b}\Big|w_k^{p+1}-w_0^{p+1}\Big|dx\notag\\
\leq&\Big(\int_{B_R}|x|^{-br}dx\Big)^{\f{1}{r}}
\Big(\int_{B_R}\Big|w_k^{p+1}-w_0^{p+1}\Big|^tdx\Big)^{\f{1}{t}}\\
\leq&C\Big(\int_{B_R}\Big|w_k^{p+1}-w_0^{p+1}\Big|^tdx\Big)^{\f{1}{t}},\notag
\end{align}
where $\f{1}{r}+\f{1}{t}=1$, $t>1$ and $r>1$ satisfies $\frac{b}{N}<\frac{1}{r}$. Note that $\f{1}{t}=1-\f{1}{r}<\frac{N-b}{N}$. Consider $p_1>0$ and $q_1>0$ satisfying
 \begin{equation}\label{eqS3.18B}
\frac{p}{p_1}+\frac{1}{q_1}=\f{1}{t}<\frac{N-b}{N},
 \end{equation}
which then yields from (\ref{eqS3.18A}) that
\begin{equation}	\label{eqS3.18C}
\arraycolsep=1.5pt\begin{array}{lll}
A_k
&\leq&C\Big(\int_{B_R}\big|w_k^{p+1}-w_0^{p+1}\big|^tdx\Big)^{\f{1}{t}}   \\[3mm]
&\le &C\big(\|w_k\|^p_{L^{p_1}(B_R)}+\|w_0\|^p_{L^{p_1}(B_R)}\big)\|w_k-w_0\|_{L^{q_1}(B_R)}.
\end{array}
\end{equation}
Similar to \cite[Theorem 1.5]{AD}, choose suitable constants $p_1>0$ and $q_1>0$ satisfying (\ref{eqS3.18B}), so that
\[
\|w_k\|^p_{L^{p_1}(B_R)}+\|w_0\|^p_{L^{p_1}(B_R)}\le C\ \ \mbox{and}\ \ \|w_k-w_0\|_{L^{q_1}(B_R)}\to 0\ \ \text{as}\ \ k\to\infty,
\]
where $C>0$ is independent of $k>0$. This further implies from (\ref{eqS3.18C}) that
\begin{equation}	\label{eqS3.18D}
A_k\to 0\ \ \text{as}\ \ k\to\infty,
\end{equation}
On the other hand, for any $\epsilon>0$, there exists $R\geq\epsilon ^{-\f{1}{b}}$ such that
\begin{align*}	
B_k&=\int_{\R^N\setminus
B_R}|x|^{-b}\Big|w_k^{p+1}-w_0^{p+1}\Big|dx\\
&\leq\epsilon\int_{\R^N\setminus
B_R}\big(w_k^{p+1}+w_0^{p+1}\big)dx\leq C\epsilon\ \ \text{as}\ \ k\to\infty,
\end{align*}
due to Sobolev's embedding theorem and the uniform boundedness of $w_k$ in $H^1(\R^N)$. Since $\epsilon>0$ is arbitrary, we conclude from above that the claim (\ref{eqS3.18}) holds true. Following \eqref{eqS3.8} and (\ref{eqS3.18}), one can deduce that $w_0\not\equiv0$.

Next, we prove that $\|w_0\|^2_2=1$. By contradiction, we assume that $\|w_0\|^2_2=l$, where $l\in(0,1)$. Set $w_l:=\f{w_0}{\sqrt{l}}$. By \eqref{eqS3.13}, we may assume that $w_k\to w_0$ $a.e.$ in $\R^N$ as $k\to\infty$. Using the Br\'{e}zis-Lieb  lemma, we obtain that
    \begin{equation}\label{11.18-1}
    \big\|\nabla w_k\big\|_2^2=\big\|\nabla  w_0 \big\|_2^2+\big\|\nabla (w_k-w_0)\big\|_2^2+o(1) \,\ \mbox{as} \,\ k\to\infty.
    \end{equation}
From \eqref{noVenergy.esti}, \eqref{energy.esti}, \eqref{eqS3.6}, (\ref{eqS3.18}) and \eqref{11.18-1}, we derive that as $k\to\infty$,
     \begin{align}\label{eqS3.15-3}
     \begin{split}
      -\lambda_0=&\lim_{k\to\infty}\epsilon^2_kI(M_k)\\
      =&\lim_{k\to\infty}\epsilon_k^2\Big(\int_{\R^N}|\nabla u_k|^2dx+\int_{\R^N}V(x)u_k^2dx-\f{2M_k^{\f{p-1}{2}}}{p+1}\int_{\R^N}\f{|u_k|^{p+1}}{|x|^b}dx\Big)\\
     \geq&\lim_{k\to\infty}\Big(\int_{\R^N}|\nabla{w}_k|^2dx
     -\f{2(a^*)^{\f{p-1}{2}}}{p+1}\int_{\R^N}\f{|{w}_k|^{p+1}}{|x|^b}dx
     \Big)\\
     =&\int_{\R^N}|\nabla w_0|^2dx+\lim_{k\to\infty}\int_{\R^N}|\nabla(w_k-w_0)|^2dx
     -\f{2(a^*)^{\f{p-1}{2}}}{p+1}\int_{\R^N}\f{|{w}_0|^{p+1}}{|x|^b}dx\\
     \geq&\int_{\R^N}|\nabla w_0|^2dx
     -\f{2(a^*)^{\f{p-1}{2}}}{p+1}\int_{\R^N}\f{|{w}_0|^{p+1}}{|x|^b}dx\\
     =&l\Big[\int_{\R^N}|\nabla w_l|^2dx-\f{2(a^*l)^{\f{p-1}{2}}}{p+1}\int_{\R^N}\f{|{w}_0|^{p+1}}{|x|^b}dx\Big]\\
     >&l \tilde{I}(a^*)=-l\lambda_0<0,
\end{split}
     \end{align}
which is a contradiction. Hence, $\|w_0\|^2_2=1$ holds true.

Since $\|w_k\|^2_2=\|w_0\|^2_2=1$, we have
\begin{equation}\label{11.18-2}
  w_k(x)\to w_0(x)\,\ \text{ strongly in}\,\ L^2(\R^N)\,\ \text{as}\,\ k\to\infty.
\end{equation}
By the weak lower semicontinuity, \eqref{noVenergy.esti} and \eqref{eqS3.18}, we then derive from \eqref{eqS3.9} that
\begin{equation}\label{eqS3.25}
     \nabla w_k(x)\to\nabla w_0(x)\, \ \text{strongly in} \,\  L^2 (\R^N)\, \ \text{as}\, \ k\to\infty.
\end{equation}
Note from \eqref{eqS3.9} that $\{w_k\}$ is a minimizing sequence of $\tilde{I}(a^*)$. One then deduces from \eqref{eqS3.18} and \eqref{eqS3.25} that $w_0$ is a minimizer of $\tilde{I}(a^*)$. By \eqref{noVminimizer}, we obtain that $w_0(x)=\f{w(x)}{\sqrt{a^*}}$.
Combining (\ref{11.18-2}) and (\ref{eqS3.25}), we obtain  that
\begin{equation}\label{11.18-3}
      w_k(x)\to w_0(x)=\f{w(x)}{\sqrt{a^*}}\, \ \text{strongly in} \,\  H^1 (\R^N)\, \ \text{as}\, \ k\to\infty,
\end{equation}
which gives \eqref{eqS3.7}. The lemma is thus proved.
\qed

Applying above lemmas, we are now ready to prove Theorem 1.1. \vskip 0.05truein

\noindent\textbf{Proof of Theorem 1.1}: 1. We first prove the exponential decay (\ref{eqSA2}). Let $u_k>0$ be a minimizer of $I(M_k)$, and consider the sequence $\{w_k\}$ defined in  Lemma \ref{lemma3.2}, where $M_k\to\infty$ as $k\to\infty$.
We claim that there exists a subsequence, still denoted by $\{w_k\}$, of $\{w_k\}$ such that
\begin{equation}\label{eqSA4}
  w_k(x)\to0\, \ \text{as}\, \ |x|\to\infty\, \ \text{uniformly for sufficiently large}\, \ k>0.
\end{equation}
Indeed, one can derive from (\ref{eqS3.7}) that for any $2\leq\alpha<2^*$,
\begin{align}\label{eqSA5}
\int_{|x|\geq\gamma}|w_k|^\alpha dx\to0\quad \text{as }\ \gamma\to\infty \ \text{uniformly for sufficiently large} \ k>0.
\end{align}
On the other hand, it follows from (\ref{A8ELE}) and (\ref{eqS3.6}) that $w_k$ satisfies the following equation
 \begin{equation}\label{eqS3.14}
  -\Delta
  w_k+\epsilon_k^2V(\epsilon_kx)w_k-(a^*)^{\f{p-1}{2}}\f{w^p_k}{|x|^b}=\mu_k\epsilon_k^2w_k \ \ \text{in}\, \ \R^N,
 \end{equation}
where ${\mu}_k\in\R$ is the Lagrange multiplier. Applying \eqref{A7wdeidentity},  (\ref{eqS3.6}) and \eqref{eqS3.18}, we deduce from Lemma \ref{Lemmaenergy.esti} that
\begin{equation}\label{eqS3.15}
\begin{split}
\epsilon_k^2\mu_k&= \epsilon_k^2\Big(I(M_k)-\f{p-1}{p+1}M_k^{\f{p-1}{2}}\int_{\R^N}\f{|u_k|^{p+1}}{|x|^b}dx\Big)\\
&=\epsilon_k^2I(M_k)-\f{p-1}{p+1}(a^*)^{\f{p-1}{2}}\int_{\R^N}\f{|{w}_k|^{p+1}}{|x|^b}dx\\
&\to-1\,\ \text{as}\,\ k\to\infty.
\end{split}
\end{equation}
Using \eqref{eqS3.15}, we derive from \eqref{eqS3.14} that as $k\to\infty$,
\begin{equation}\label{11.19-1}
-\Delta w_k-c(x)w_k\leq0\ \ \mbox{in}\,\ \R^N,\,\ \text{where}\  \ c(x)=(a^*)^{\f{p-1}{2}}\f{w^{p-1}_k(x)}{|x|^b}.
\end{equation}
Furthermore, one can check from H\"{o}lder inequality that
$$c(x)\in L^t(\R^N), \ \ \mbox{where}\ \, t\in\Big(\f{2N}{N(p-1)+2b},\f{2N}{(N-2)(p-1)+2b}\Big).$$
Applying De Giorgi-Nash-Moser theory (\cite[Theorem 4.1]{QF}) to \eqref{11.19-1}, we deduce that
\begin{align}\label{eqSA6}
\underset{B_1(\xi)}{\text{max}}w_k(x)\leq C\Big(\int_{B_2(\xi)}|w_k(x)|^\alpha dx\Big)^\f{1}{\alpha}\quad \text{for sufficiently large}\, \ k>0,
\end{align}
where $\xi\in\R^N$ is arbitrary, and $C>0$ depends only on the bound of $\|c(x)\|_{L^t(B_2(\xi))}$. Thus, (\ref{eqSA4}) follows from (\ref{eqSA5}) and (\ref{eqSA6}).

Due to the smallness of $|x|^{-b}$ for large $|x|>0$, we now derive from (\ref{eqSA4}), \eqref{eqS3.14} and \eqref{eqS3.15} that there exists a sufficiently large constant $R>0$, which is independent of $k$, such that as $k\to\infty$,
\begin{equation}\label{eqSA8}
  -\Delta w_k(x)+\theta w_k(x)\leq0 \ \ \text{in}\ \ \R^N\setminus B_R(0),
 \end{equation}
where $0<\theta<1$ is independent of $k$. By the comparison principle \cite [Theorem 6.4.2]{CPY}, we obtain from (\ref{eqSA8}) that as $k\to\infty$,
\begin{equation}\label{eqSA9}
w_k(x)\leq Ce^{-\sqrt{\theta}|x|}\ \ \, \text{in}\ \ \R^N\setminus B_R(0),
 \end{equation}
which implies the exponential decay (\ref{eqSA2}) for $w_k$ as $k\to\infty$. Moreover, under the assumption $(V)$, since the term $|x|^{-b}$ is small for large $|x|$, applying the local elliptic estimate (cf. (3.15) in \cite{GT}) yields from (\ref{eqSA9}) that as $k\to\infty$,
$$|\nabla w_k(x)|\leq Ce^{-\theta|x|}\ \ \, \mbox{for}\ \ |x|>R,$$
which thus gives the exponential decay (\ref{eqSA2}) for $\nabla w_k$ as $k\to\infty$. This proves (\ref{eqSA2}).

2. We next prove that (\ref{A9limit.behavior}) holds true. Recall from Lemma \ref{lemma3.2} that
\begin{equation}\label{eqSA10}
w_k(x)\to \f{w(x)}{\sqrt{a^*}} \,\ \text{in}\, \ H^1(\R^N)\, \ \text{as} \,\ k\to\infty,
\end{equation}
where the convergence holds for the whole sequence $\{w_k(x)\}$, due to the uniqueness of $w(x)>0$. Following (\ref{eqSA10}), the $L^\infty$-uniform convergence (\ref{A9limit.behavior}) for the case $N=1$ can be directly obtained by applying Sobolev's embedding theorem $H^1(\R)\hookrightarrow L^\infty(\R)$.

We now prove the $L^\infty$-uniform convergence (\ref{A9limit.behavior}) for the case $N\geq2$. 
Rewrite \eqref{eqS3.14} as
 \begin{equation}\label{T1}
  -\Delta
  w_k(x)=G_k(x) \ \ \text{in}\, \ H^1(\R^N),
 \end{equation}
where $$G_k(x):=\mu_k\epsilon_k^2w_k-\epsilon_k^2V(\epsilon_kx)w_k+(a^*)^{\f{p-1}{2}}\f{w^p_k}{|x|^b}.$$
Since $w_k$ is bounded uniformly in $H^1(\R^N)$ as $k\to\infty$, we derive from \eqref{eqSA6} that
 \begin{equation}\label{T2}
  w_k \, \ \text{is bounded uniformly in}\, \ L^\infty(\R^N).
 \end{equation}
Using H\"{o}lder inequality, we deduce that for any $R>0$,
\begin{equation}\label{eqSA12}
\int_{B_R(0)}\Big|\f{w^p_k(x)}{|x|^b}\Big|^rdx
\leq\Big(\int_{B_R(0)}\f{1}{|x|^{brt}}dx\Big)^{\f{1}{t}}
\Big(\int_{B_R(0)}|w^p_k|^{rt'}dx\Big)^{\f{1}{t'}}\\
<\infty,
\end{equation}
where $r\in(1,\f{N}{b})$, $\f{1}{t}+\f{1}{t'}=1$, and $t'=1+\max\Big\{\f{N}{N-br},\f{2}{pr}\Big\}$. We then obtain from \eqref{eqSA12} that
\begin{equation}\label{eqSA14}
w^p_k(x)|x|^{-b}\in L^r_{loc}(\R^N)\quad \text{for any}\ r\in \big(1,\f{N}{b}\big),
\end{equation}
which implies that $G_k(x)$ is bounded uniformly in $L^r_{loc}(\R^N)$. For any large $R>0$, it thus follows from \cite [Theorem 9.11]{GT} that
\begin{equation}
  \|w_k(x)\|_{W^{2,r}(B_R)}\leq C\Big(\|w_k(x)\|_{L^r(B_{R+1})}+\|G_k(x)\|_{L^r(B_{R+1})}\Big),
\end{equation}
where $C>0$ is independent of $k>0$ and $R>0$. By the compactness of the embedding
$W^{2,r}(B_R)\hookrightarrow L^\infty(B_R)$ for $2r>N$, cf. \cite [Theorem 7.26] {GT}, we conclude that there exists a subsequence, still denoted by $\{w_k\}$, of $\{w_k\}$ such that
\begin{equation}\label{T3}
w_k(x)\to \tilde{w}_0(x)\, \ \text{uniformly in}\, \ L^\infty(B_R)\, \ \text{as} \,\ k\to\infty.
\end{equation}
Since $R>0$ is arbitrary, we obtain from \eqref{eqSA10} that
\begin{equation}\label{T3}
w_k(x)\to \f{w(x)}{\sqrt{a^*}}\, \ \text{uniformly in}\, \ L^\infty_{loc}(\R^N)\, \ \text{as} \,\ k\to\infty.
\end{equation}

On the other hand, we deduce from (\ref{A6wdedecay}) and (\ref{eqSA2}) that for any $\epsilon>0$, there exists a constant
$R_\epsilon>0$, independent of $k>0$, such that
$$|w_k(x)|, \ \,\big|\f{w(x)}{\sqrt{a^*}}\big|<\f{\epsilon}{4} \quad \text{for any} \ |x|>R_\epsilon,$$
which implies that
$$\sup_{|x|>R_\epsilon}\big|w_k(x)-\f{w(x)}{\sqrt{a^*}}\big|\leq
\sup_{|x|>R_\epsilon}\big(|w_k(x)|+\big|\f{w(x)}{\sqrt{a^*}}\big|\big)
\leq\f{\epsilon}{2}.$$
Recall from (\ref{T3}) that for sufficiently large $k>0$,
$$\sup_{|x|\leq R_\epsilon}\big|w_k(x)-\f{w(x)}{\sqrt{a^*}}\big|
\leq\f{\epsilon}{2}.$$
We now conclude from above that the $L^\infty$-uniform convergence (\ref{A9limit.behavior}) holds true for all $N\geq2$. The proof of Theorem \ref{limit.theorem} is therefore complete.
\qed

\vspace{.9cm}
%


\end{document}